\documentclass{article}
 
\usepackage{amsmath}
\usepackage{amssymb}
\usepackage{amsthm}
\usepackage[usenames,dvipsnames,svgnames,table]{xcolor}
\usepackage[colorlinks,linktocpage, backref, colorlinks=true,
            linkcolor=red,
            urlcolor=blue,
            citecolor=gray]{hyperref}
\usepackage{amsrefs}

\DeclareMathOperator{\homdist}{HomDist}
\DeclareMathOperator{\dist}{dist}

\DeclareMathOperator{\Hom}{Hom}
\DeclareMathOperator{\defect}{def}
\DeclareMathOperator{\rk}{rk}
\DeclareMathOperator{\PSL}{PSL}
\DeclareMathOperator{\SO}{SO}
\DeclareMathOperator{\SU}{SU}
\DeclareMathOperator{\PU}{PU}

\def\Z{\mathbb Z}
\def\N{\mathbb{N}}
\def\F{{\mathbf F}}
\newcommand{\Uni}{\mathrm{U}}
\newcommand{\uf}{\mathcal{U}}
\newcommand{\Uniuf}{\Uni_\uf}
\newcommand{\Matuf}{\mathrm{M}_\uf}
\newcommand{\freegrp}{\mathbf{F}}
\newcommand{\subgrpeq}{\leq}
\newcommand{\norsubgrpeq}{\trianglelefteq}
\newcommand{\commutator}[1]{[#1]}
\newcommand{\catFin}{\mathbf{Fin}}
\newcommand{\catAlt}{\mathbf{Alt}}
\newcommand{\catNil}{\mathbf{Nil}}
\newcommand{\catSol}{\mathbf{Sol}}

\newcommand{\catPSL}{\mathbf{PSL}}
\newcommand{\reals}{\mathbb{R}}
\newcommand{\calC}{\mathcal{C}}
\newcommand{\calU}{\mathcal{U}}

\newcommand{\frob}[1]{\lVert #1 \rVert_{\mathrm{Frob}}}
\newcommand{\opnorm}[1]{\lVert #1 \rVert_{\mathrm{op}}}
\newcommand{\hsnorm}[1]{\lVert #1 \rVert_{\mathrm{HS}}}
\newcommand{\closure}[1]{\overline{#1}}
\newcommand{\card}[1]{\lvert#1\rvert}

\theoremstyle{plain}
\newtheorem{lemma}{Lemma}[section]
\newtheorem{theorem}[lemma]{Theorem}
\newtheorem{definition}[lemma]{Definition}
\newtheorem{corollary}[lemma]{Corollary}
\newtheorem{proposition}[lemma]{Proposition}

\theoremstyle{definition}

\newtheorem{question}[lemma]{Question}

\title{Finitary approximations of groups  \\ and their applications}

\author{Andreas Thom \\ Institut f\"ur Geometrie, TU Dresden}

\begin{document}

\maketitle

\begin{abstract}
In these notes we will survey recent results on various finitary approximation properties of infinite groups. We will discuss various restrictions on groups that are approximated for example by finite solvable groups or finite-dimensional unitary groups with the Frobenius metric. Towards the end, we also briefly discuss various applications of those approximation properties to the understanding of the equational theory of a group. 
\end{abstract}

\setcounter{tocdepth}{2}

\section{Introduction}
\subsection{The setup}
Let $\Gamma$ be a finitely presented group, given by a finite generating set $X:=\{ x_1,\dots,x_k\}$ and a finite set $R \subset {\mathbf F}_k = \langle X \rangle$ of relations, i.e.\ $\Gamma:= \langle X | R \rangle = {\mathbf F}_k/N$, where ${\mathbf F}_k$ denotes the free group on $X$ and $N = \langle \! \langle R \rangle \! \rangle$ the normal subgroup generated by $R$. Throughout the article, $X$ and $R$ will be fixed.

Very basic questions about $\Gamma$ are usually hard to answer unless $\Gamma$ can be realised as a group of symmetries of a sufficiently concrete object, such as a finite set, a finite-dimensional vector space or a metric space with suitable properties. In the easiest situation, maybe $\Gamma$ is residually finite, i.e.\ for any finite subset $F \subset \Gamma$, there exists a homomorphism to a finite group $\varphi \colon \Gamma \to H$, such that the restriction of $\varphi$ to $F$ is injective.
In order to overcome the algebraic and arithmetic obstruction to the existence of finite quotients and finite-dimensional unitary representations it is worthwhile to relax these notions. Informally speaking, we will seek for asymptotic homomorphisms from $\Gamma$ with values in a family of (typically compact) metric groups.
Let $(G,d)$ be a metric group and assume throughout the entire article that $d \colon G \times G \to [0,\infty)$ is bi-invariant, i.e.\ 
$$d(gh,gk)=d(h,k)=d(hg,kg), \quad \forall g,h,k \in G.$$ Note that any bi-invariant metric is uniquely determined by the associated \emph{invariant length function} $\ell(g)=d(1,g)$, which is a subadditive, symmetric and conjugation invariant, $[0,\infty)$-valued function on $G$ that takes the value $0$ only at $1_G$.

Well-known  invariant length functions in this context include the normalized Hamming distance on symmetric groups or various length functions induced by unitarily invariant norms on groups of unitary matrices. Somewhat similarly, there is a {rank length function}, which is defined if $G\subgrpeq {\rm GL}_n(q)$ for some $n\in\N$ and $q$ a prime power:
$
\ell_G^{\rm rank}(g):= n^{-1} \cdot {\rk(1-{g})}.
$
Another important example is the {conjugacy length function}, which is defined by
\begin{equation*} \label{conjmetric}
\ell^{\rm conj}(g):= \log_{|G|}\card{g^G}
\end{equation*}
for $g\in G$, where $g^G$ is the conjugacy class of $g \in G$, and $G$ is a finite and centerless group. If $G$ is the alternating group or a finite simple group of Lie type, then the conjugacy metric is comparable to the more geometrically defined metrics above. A fundamental result in the work of Liebeck-Shalev \cite{liebeckshalev2001diameters} says 
that the the conjugacy length is intrinsically tied to the algebraic properties of $G$ if $G$ is simple: indeed, there exists an absolut constant $c>0$, such that $(g^G)^{ck} = G$ if $k > \ell^{\rm conj}(G)^{-1}$. Thus, the conjugacy length is also comparable to the normalized word metric w.r.t.\ any sufficiently small conjugacy class. It turns out that there is essentially just one invariant length function up to a suitable notion of equivalence on a finite simple group.

\subsection{Approximation and stability }
Now we can define more precisely what we mean by metric approximation of an abstract group by a class of metric groups $\calC$. 

\begin{definition}\label{def:C_approx_abs_grp}
A group $\Gamma$ is called \emph{$\calC$-approximated} if there is a length function $\delta\colon\Gamma\to[0,\infty)$ such that for any finite subset $S\subseteq \Gamma$ and $\varepsilon>0$ there exist a group $(G,d)\in\calC$ and a map $\varphi \colon \Gamma \to G$, such that
	\begin{enumerate}
		\item[$(i)$] if $g,h,gh\in S$, then $d(\varphi(g)\varphi(h),\varphi(gh))<\varepsilon$ and
		\item[$(ii)$] for $g\in S$ we have $d(1_H,\varphi(g))\geq\delta(g)$.
	\end{enumerate}
\end{definition}

In some situations, we will only fix a class $\calC$ of groups and let the choice of bi-invariant metrics be arbitrary -- we will also speak about $\calC$-approximability in this context. The previous definition has emerged from various contexts, including an influential work of Ulam \cite{ulam}, the work of Connes \cite{MR0454659}, Gromov \cite{gromov1999endomorphisms}, Weiss \cite{weiss } and later work \cite{MR2900231, holtrees2016some, glebsky2016approximations}.

Let us also introduce the closely related notion of asymptotic homomorphisms.
Note that any map $\varphi\colon X \to G$, for some $(G,d)\in \calC$, uniquely determines a homomorphism ${\mathbf F}_k\to G$ which we will also denote by $\varphi$.

\begin{definition}
Let $(G,d)\in\calC$ and let $\varphi,\psi\colon X\to G$ be maps. The \emph{defect} of $\varphi$ is defined by
\[\defect(\varphi):=\max_{r\in R}d(\varphi(r),1_G).\]
The \emph{distance} between $\varphi$ and $\psi$ is defined by \[\dist(\varphi,\psi)=\max_{1 \leq i \leq k}d(\varphi(x_i),\psi(x_i)).\]
The \emph{homomorphism distance} of $\varphi$ is defined by
\[\homdist(\varphi):=\inf_{\pi\in\Hom(\Gamma,G)}\dist(\varphi,\pi|_{X}).\]
\end{definition}

\begin{definition}
A sequence of maps $\varphi_n\colon  X\to G_n$, for $(G_n,d_n)\in\calC$, is called an \emph{asymptotic homomorphism} with values in $\calC$ if 
$$\lim_{n\to \infty}\defect(\varphi_n)=0.$$
\end{definition}

\begin{definition}
Let $G_n\in\calC,n\in \N$. Two sequences $\varphi_n,\psi_n\colon  X\to G_n$ are called equivalent if $$\lim_{n\to\infty}\dist(\varphi_n,\psi_n)=0.$$ If an asymptotic homomorphism $(\varphi_n)_{n\in\N}$ is equivalent to a sequence of homomorphisms, we call $(\varphi_n)_{n\in\N}$ \emph{trivial}.
\end{definition}

It is easy to see that a group $\Gamma= \langle X | R \rangle$ is $\calC$-approximated if and only if there is an asymptotic homomorphism with values in $\calC$ that separates elements in a suitable sense. Note that the existence of a finite presentation is assumed mostly for convenience. Moreover, note that it is easy to see that the property of being $\calC$-approximated depends only on $\Gamma$ and not on the finite presentation.

Let us discuss some examples of $\calC$-approximated groups. We denote by $\catAlt$ (resp.~$\catFin$) the class of finite alternating groups (resp.~the class of all finite groups). 
A group is called sofic (resp.~weakly sofic) if and only if it is $\catAlt$-approximated (resp.~$\catFin$-approximated) as an abstract group, see \cite{glebsky2016approximations}.
The class of sofic groups is of central interest in group theory. Indeed, eversince the work of Gromov on Gottschalk's Surjunctivity Conjecture \cite{gromov1999endomorphisms}, the class of sofic groups has attracted much interest in various areas of mathematics. Major applications of this notion arose in the work Elek and Szab\'o  on Kaplansky's Direct Finiteness Conjecture \cite{elekszabo2004sofic}, L\"uck's Determinant Conjecture \cite{elekszabo2005hyperlinearity}, and more recently in joint work of the author with Klyachko on generalizations of the Kervaire-Laudenbach Conjecture and Howie's Conjecture \cite{MR3604379}. 
Despite considerable effort, no non-sofic group has been found so far -- whether all groups are sofic is one of the outstanding open problems in group theory. 

\begin{question}[Gromov]
Are all groups sofic?
\end{question}

Examples of sofic groups which fail to be locally residually amenable are given in \cite{cornulier2011sofic} and \cite{karnikolov2014non} (see also \cite{MR2566306}), answering a question of Gromov \cite{gromov1999endomorphisms}.

Groups approximated by certain classes of finite simple groups of Lie type have been studied in \cite{arzhantsevapaunescu2017linear} and \cite{MR3210125, 1606.03863}. We will discuss approximation by groups of unitary matrices, a central topic in the theory of operator algebras and free probability theory, at length in Section \ref{unitary}.

Non-approximation results are rare, however, in \cite{MR2900231} it was proved that the so-called Higman group cannot be approximated by finite groups with commutator-contractive invariant length functions. In \cite{howie1984the} Howie presented a group which (by a result of Glebsky \cite{glebskyrivera2008sofic}) turned out not to be approximated by finite nilpotent groups with arbitrary invariant length function. In Sections \ref{finitegrp} and \ref{unitary}, we will survey more general results of this type that have been proved recently in \cite{1703.06092}.

A central definition in the present context is the notion of stability that was introduced in \cite{lubchigletho}.

\begin{definition}
The group $\Gamma$ is called \emph{$\calC$-stable} if all asymptotic homomorphisms with values in $\calC$ are trivial, that is: for all $\varepsilon>0$ there exists $\delta>0$ such that
$\defect(\varphi)< \delta$ implies $\homdist(\varphi)< \varepsilon$ for all $\varphi \colon X \to (G,d)$ with $(G,d) \in \calC$.
\end{definition}

It is clear that any group $\Gamma$ that is both $\calC$-approximated and $\calC$-stable is residually a subgroup of groups in $\calC$. If $\calC$ consists of finite groups or more generally compact groups this readily implies that $\Gamma$ must be residually finite. This observation has been used in \cite{lubchigletho} to prove non-approximation results by proving that certain groups are $\calC$-stable but not residually $\calC$, see Theorem \ref{main}.

\subsection{Metric ultraproducts}\label{ultra}

Throughout these notes, we fix a non-principal ultrafilter $\uf$ on $\N$. 
Let $(G_n)_{n\in\N}$ be a family of groups, all equipped with bi-invariant metrics $d_n$.
In this case, the subgroup
\[N=\left\{(g_n)_{n\in\N}\in\prod_{n\in\N}G_n| \lim_{n\to\uf}d_n(g_n,1_{G_n})=0\right\}\]
of the direct product $\prod_{n\in\N} G_n$
is normal, so that we can define the \emph{metric ultraproduct}
\[\prod_{n\to\uf}(G_n,d_n):=\prod_{n\in\N}G_n\Big/N.\]

The relevance of metric ultraproducts becomes apparent in the following folklore result:

\begin{proposition}\label{lem:char_C_approx_abs_grp_via_ultraprod}
Let $\calC$ is a class of metric groups. A group $\Gamma = \langle X | R \rangle$ is $\calC$-approximated if and only if it is isomorphic to a subgroup of a metric ultraproduct of $\calC$-groups.
\end{proposition}

For more details on the algebraic and geometric structure of such ultraproducts see also \cite{MR3162821} and \cite{MR3210125, 1606.03863, 1709.06286}.
In view of Proposition \ref{lem:char_C_approx_abs_grp_via_ultraprod} it is natural to generalize the notion of a $\calC$-approximated group to topological groups using ultraproducts:

\begin{definition}\label{def:C_approx_top_grp}
A topological group is called \emph{$\calC$-approximated} if it is topologically isomorphic to a closed subgroup of a metric ultraproduct of $\calC$-groups.
\end{definition}

We will constrain ourselves to Polish groups and countable ultraproducts, but that is just for convenience. Typically, for example in the context of sofic or weakly sofic groups, it is easy to see that an abstract $\calC$-approximated group is also $\calC$-approximated when viewed as a topological group with the discrete topology. It is clear that any $\catFin$-approximated topological group must admit a bi-invariant metric that induces the topology. We will discuss various less obvious restrictions on $\catFin$-approximability in Section \ref{sec:approx_Lie_grps}.

In view of the definition of metric ultraproducts, any approximation property for a group $\Gamma$ by a class of compact groups leads to an embedding into a quotient of a compact group.

\begin{question}
Is any group a sub-quotient of a compact group?
\end{question}

\section{Weak soficity and the pro-finite topology} \label{finitegrp}

\subsection{Connections with the pro-finite topology}

In this section, we want to survey some recent results that were proved in joint work with Nikolov and Schneider, see \cite{1703.06092}. The main insight that helped us was to combine the relationship between soficity and properties of the pro-finite topology that was established by Glebsky-Rivera \cite{glebskyrivera2008sofic} with the deep work on finite groups by Nikolov-Segal \cite{segal2009words, MR2995181
}.
Let $\calC$ be a class of finite groups. Adapting Theorem 4.3 of \cite{glebskyrivera2008sofic}, one can prove the following theorem relating $\calC$-approximated groups to properties of the pro-$\calC$ topology on a free group:

\begin{theorem}[\cite{1703.06092}]\label{thm:char_C_approx_grps_pro_C}
	Let $\freegrp_k/N$ be a presentation of a group $\Gamma$.
	Then, if $\Gamma$ is $\calC$-approximated, for each finite sequence $n_1,\ldots,n_m\in N$ it holds that $$\closure{n_1^\freegrp\cdots n_m^\freegrp}\subseteq N,$$ where the closure is taken in the pro-$\calC$ topology on $\freegrp_k$. The converse holds under mild assumptions on $\calC$.
\end{theorem}

The coarsest such topology on $\freegrp_k$ is of course the pro-finite topology and at the time of writing of \cite{glebskyrivera2008sofic} it was an open problem to decide whether a finite product of conjugacy classes in a non-abelian free group is always closed in this topology. As has been remarked in \cite{1703.06092}, it is a rather straightforward consequence of the work of Nikolov-Segal (see \cite{MR2995181} or Theorem \ref{thm:gen_com_fin_grp}) that this is not the case. Indeed, one of their main results implies that in $\freegrp_k$ the profinite closure of a finite product of conjugacy classes of $x_1^{-1},x_1,\ldots,x_k^{-1},x_k$ contains the entire commutator subgroup, while it is a well known fact (see Theorem 3.1.2 of \cite{segal2009words}) that the commutator width of $\freegrp_k$ is infinite if $k>1$. 
This implication was first observed by Segal and independently discovered by Gismatullin. In view of this observation it seems unlikely that the pro-finite closure is always contained in the normal closure, but this remains an open problem.

\begin{question}[Glebsky-Rivera, \cite{glebskyrivera2008sofic}] \label{q1}
Let $n_1,\dots,n_m \in \freegrp_k$. Is it true, that
$$\closure{n_1^\freegrp\cdots n_m^\freegrp}\subseteq  \langle \! \langle n_1,\dots,n_m \rangle \! \rangle,$$ where the closure is taken in the pro-finite topology on $\freegrp_k$?
\end{question}

In general, there are quite a number of mysteries that can be formulated in terms of closure properties of the pro-$\calC$ topology for particular more restricted families of groups. Indeed, let us just mention a question from \cite{MR1621745}. 

\begin{question}[Herwig-Lascar, \cite{MR1621745}]
It is easy to see that if a finitely generated subgroup $H < \freegrp_k$ is closed in the pro-odd topology, then it satisfies $a^2 \in H \Rightarrow a \in H$. Is the converse true?
\end{question}

\subsection{Approximation by classes of finite groups}

For more restricted families of groups, the answer to Question \ref{q1} becomes negative. Indeed,
let $\catSol$ (resp.~$\catNil$) be the class of finite solvable (resp.~nilpotent) groups. In view of Theorem \ref{thm:char_C_approx_grps_pro_C} this implies that there are groups which are not $\catSol$-approximated. More precisely, we proved:

\begin{theorem}[\cite{1703.06092}]\label{thm:fin_gen_perf_grp_not_C_approx}
Every finitely generated $\catSol$-approximated group has a non-trivial abelian quotient.
\end{theorem}

As a consequence, perfect groups cannot be $\catSol$-approximated and a finite group is $\catSol$-approximated if and only if it is solvable. Indeed, any finite solvable group is clearly $\catSol$-approximated and on the other hand, a non-solvable finite group contains a non-trivial perfect subgroup and hence cannot be $\catSol$-approximated by  Theorem \ref{thm:fin_gen_perf_grp_not_C_approx}.

Initially, Howie proved in \cite{howie1984the} that the group $\langle  x,y | x^{-2}y^{-3},x^{-2}(xy)^5 \rangle$ is not $\catNil$-approximated. We followed his proof for any non-trivial finitely generated perfect group and then extended it in \cite{1703.06092} and established that these groups are not even $\catSol$-approximated using techniques of Segal \cite{segal2000closed,segal2009words}.

Note that finite generation is crucial in the statement of Theorem \ref{thm:fin_gen_perf_grp_not_C_approx}. Indeed, there exist countably infinite locally finite-$p$ groups which are perfect and even characteristically simple, see \cite{mclain1954characteristically}. These groups are $\catNil$-approximated, since finite $p$-groups are nilpotent. It is known that locally finite-solvable groups cannot be non-abelian simple, but it seems to be an open problem if there exist $\catSol$-approximated simple groups.

Let us also remark, that the assumptions of Theorem \ref{thm:fin_gen_perf_grp_not_C_approx} are not enough to conclude that a finitely generated $\catSol$-approximated group has an infinite solvable quotient. Indeed, consider a suitable congruence subgroup of ${\rm SL}(3,\Z)$, which is residually $p$-finite and thus even $\catNil$-approximated and has Kazhdan's property (T) -- thus all amenable quotients are finite. However, the following seems to be an open problem.

\begin{question}
Is every finitely presented and $\catSol$-approximated group residually finite-solvable?
\end{question}

Even more, it could be that all finitely presented groups are $\catSol$-stable (in a suitable sense). A positive answer to the previous question would be in sharp contrast to other forms of approximability. Indeed, for example the Baumslag-Solitar group ${\rm BS}(2,3)$ is not residually finite, but residually solvable and hence sofic and in particular $\catFin$-approximated.

Let us finish this section by mentioning some structure result on the class of $\catFin$-approximated groups. Let $\catPSL$ be the class of simple groups of type $\PSL(n,q)$, i.e.~$n\in\N_{\geq 2}$ and $q$ is a prime power and $(n,q)\neq (2,2),(2,3)$, and recall that $\catFin$ is the class of all finite groups.
In \cite{1703.06092} we prove the following result.

\begin{theorem}[\cite{1703.06092}] \label{thm:smpl_Fin_approx_grp_ultraprod_smpl_grps}
	Any non-trivial finitely generated $\catFin$-approximated group has a non-trivial $\catPSL$-approximated quotient. In particular, every finitely generated simple and $\catFin$-approximated group is $\catPSL$-approximated.
\end{theorem}

The proof makes use of seminal results of Liebeck-Shalev \cite{liebeckshalev2001diameters} and Nikolov-Segal \cite{MR2995181}. The previous result may be seen as a first step towards a proof that all $\catFin$-approximated groups are sofic.

\subsection{Approximability of Lie groups}\label{sec:approx_Lie_grps}

Let us explain how a theorem of Nikolov-Segal allows us to deduce two results concerning the approximability of Lie groups by finite groups and one result on compactifications of pseudofinite groups.

\begin{theorem}[Theorem 1.2 of \cite{MR2995181}]\label{thm:gen_com_fin_grp}
	Let $g_1,\ldots,g_m$ be a symmetric generating set for the finite group $G$. If $K\norsubgrpeq G$, then
	$$
	\commutator{K,G}={\left(\prod_{j=1}^m{\commutator{K,g_j}}\right)}^{e},
	$$
	where $e$ only depends on $m$.
\end{theorem}

The following result is an immediate corollary of Theorem \ref{thm:gen_com_fin_grp}.

\begin{corollary}[\cite{1703.06092}]\label{cor:two_gen_com_fin_grp}
	Let $G$ be a finite group, then for $g,h\in G$ and $k\in\N$ we have
	$$
	\commutator{g^k,h^k}\in {\left(\commutator{G,g}\commutator{G,g^{-1}}\commutator{G,h}\commutator{G,h^{-1}}\right)}^{e}
	$$
	for some fixed constant $e\in\N$ that is independent of $G$.
\end{corollary}

We deduce immediately that the same conclusion holds for any quotient of a product of finite groups and in particular, for any metric ultraproduct of finite groups. Combining the finitary approximation with the local geometry of Lie groups we obtain the following consequence.

\begin{theorem}[\cite{1703.06092}]\label{thm:ctd_Fin_approx_Lie_grps_ab}
	A connected Lie group is $\catFin$-approximated as a topological group if and only if it is abelian.
\end{theorem}

Indeed, this is a direct consequence of the following auxiliary result:

\begin{lemma}\label{lem:two_one_par_subgrps_to_ultprod_comm}
	Let $\varphi,\psi:\reals\to(H,\ell)=\prod_{\calU}{(H_i,\ell_i)}$ be continuous homomorphisms into a metric ultraproduct of finite metric groups with bi-invariant metrics. Then, the images of $\varphi$ and $\psi$ commute.
\end{lemma}

Note that Theorem \ref{thm:ctd_Fin_approx_Lie_grps_ab} provides an answer to Question 2.11 of Doucha \cite{doucha2016metric} whether there are groups with invariant length function that do not embed in a metric ultraproduct of finite groups with invariant length function. Note also that the topology matters a lot in this context. Indeed, it can be shown that any compact Lie group is a discrete subgroup of a countable metric ultraproduct of finite groups, see \cite{1703.06092}.

When one restricts the class of finite groups further and approximates with symmetric groups, one can not even map the real line $\reals$ non-trivially and continuously to a metric ultraproduct of such groups with invariant length function. Indeed, for the symmetric group ${\rm Sym}(n)$, it can be shown that all invariant length functions $\ell$ on it satisfy $\ell(\sigma^k)\leq 3\ell(\sigma)$, for every $k\in\Z$ and $\sigma\in {\rm Sym}(n)$. Using this identity, it is simple to deduce that the only continuous homomorphism of $\reals$ into a metric ultraproduct of finite symmetric groups with invariant length function is trivial.

Referring to a question of Zilber \cite[p.~17]{zilber2014perfect} (see also Question 1.1 of Pillay \cite{pillay2015remarks}) whether a compact simple Lie group can be a quotient of the algebraic ultraproduct of finite groups, we obtained the following second application of Corollary \ref{cor:two_gen_com_fin_grp}:

\begin{theorem}[\cite{1703.06092}]\label{thm:Lie_grp_quot_prod_fin_grps}
Let $G$ be a Lie group equipped with an bi-invariant metric generating its topology. If $G$ is an abstract quotient of a product of finite groups, then $G$ has abelian identity component.
\end{theorem}

The proof of this result is almost identical to the proof of Theorem \ref{thm:ctd_Fin_approx_Lie_grps_ab}.
Theorem \ref{thm:Lie_grp_quot_prod_fin_grps} implies that any compact simple Lie group, the simplest example being $\SO(3,\reals)$, is not a quotient of a product of finite groups, answering the questions of Zilber and Pillay. Note also that these results are vast generalizations of an ancient result of Turing {\cite{turing1938finite}.

Moreover, Theorem \ref{thm:Lie_grp_quot_prod_fin_grps} remains valid if we replace the product of finite groups by a \emph{pseudofinite group}, i.e.~a group which is a model of the theory of all finite groups.
It then also provides a negative answer to Question 1.2 of Pillay \cite{pillay2015remarks}, whether there is a surjective homomorphism from a pseudofinite group to a compact simple Lie group.

By a \emph{compactification} of an abstract group $G$, we mean a compact group $C$ together with a homomorphism $\iota \colon G\to C$ with dense image.
Pillay conjectured that the Bohr compactification (i.e.~the universal compactification) of a pseudofinite group has abelian identity component (Conjecture 1.7 in \cite{pillay2015remarks}).
We answer this conjecture in the affirmative by the following result:
\begin{theorem}[\cite{1703.06092}]\label{thm:cpt_pseudo_fin_grp_ab}
	Let $G$ be a pseudofinite group. Then the identity component of any compactification $\iota \colon G\to C$ is abelian.
\end{theorem}
Again, the proof is an application of Corollary \ref{cor:two_gen_com_fin_grp}.

\section{Approximation by unitary matrices} \label{unitary}

\subsection{The choice of the metric}

We will now focus on approximation of groups by unitary matrices.
Today, the theme knows many variations, ranging from operator-norm approximations that appeared in the theory of operator algebras \cite{MR1437044, cde} to questions related to Connes' Embedding Problem, see \cite{MR0454659, MR2460675} for details. 
Several examples of this situation have been studied in the literature:
\begin{enumerate}
\item[(1)] $G_n={\rm U}(n)$, where the metric $d_n$ is induced by the Hilbert-Schmidt norm $\hsnorm T = \sqrt{n^{-1} \sum_{i,j=1}^n |T_{ij}|^2}$. In this case, approximated groups are sometimes called hyperlinear \cite{MR2460675}, but we choose to call them Connes-embeddable.
\item[(2)] $G_n={\rm U}(n)$, where the metric $d_n$ is induced by the operator norm $\|T\|_{\rm op} = \sup_{\|v\|=1} \|Tv\|.$ In this case, groups which are $(G_n,d_n)_{n=1}^{\infty}$-approximated are called {\rm MF}, see \cite{cde}.
\item[(3)] $G_n={\rm U}(n)$, where the metric $d_n$ is induced by the \emph{unnormalized} Hilbert-Schmidt norm $\|T\|_{\rm Frob}= \sqrt{ \sum_{i,j=1}^n |T_{ij}|^2}$, also called Frobenius norm. We will speak about Frobenius-approximated groups in this context, see \cite{lubchigletho}.
\end{enumerate}

Let us emphasize that the approximation properties are \emph{local} in the sense that only finitely many group elements and their relations have to be considered at a time. This is in stark contrast to the \emph{uniform} situation, which -- starting with the work of Grove-Karcher-Ruh and Kazhdan \cite{johnson, art:kazh} -- is much better understood, see \cite{MR3038548, 1706.04544}.

Again, there are longstanding problems that ask if \emph{any} group exists which is not approximated in the sense of (1), a problem closely related to Connes' Embedding Problem \cite{MR0454659, MR2460675}. 

\begin{question}[Connes, \cite{MR0454659}]
Is every discrete group Connes-embeddable?
\end{question}

Connes' Embedding Problem has many incarnations and we want to mention only a few of them, see \cite{MR2072092, MR2460675} for more details. The most striking alternative formulation is due to Kirchberg, who showed that Connes' Embedding Problem has an affirmative answer if and only if the group $\freegrp_2 \times \freegrp_2$ is {residually finite dimensional}, i.e.\ if the finite-dimensional unitary representations of this group are dense in the unitary dual equipped with the Fell topology.

Kirchberg \cite{MR1437044} conjectured that any stably finite $C^*$-algebra is embeddable into an norm-ultraproduct of matrix algebras, implying a positive answer to the approximation problem in the sense of (2) for any group. Recent breakthrough results imply that any amenable group is MF, i.e.\ approximated in the  sense of (2), see \cite{MR3583354}. 

Approximation in the sense of (3) is known to be more restrictive -- as has been shown in \cite{lubchigletho}.
Indeed, in joint work with De Chiffre,  Glebsky, and Lubotzky \cite{lubchigletho}, a conceptually new technique was introduced, that allowed to provide groups that are not approximated in the sense of (3) above. An analogous result for the normalized Frobenius norm would answer the Connes' Embedding Problem. Even though we had little to say about Connes' Embedding Problem, we believe that we provided a promising new angle of attack.

\begin{theorem}[\cite{lubchigletho}] \label{main}
There exists a finitely presented group, which is not Frobenius-approximated. Specifically, we can take a certain central extension of a lattice in ${\rm U}(2n) \cap {\rm Sp}(2n,{\mathbb Z}[i,1/p])$ 
for a large enough prime $p$ and $n \geq 3$.
\end{theorem}

The key insight was that there exists a cohomological obstruction (in the second cohomology with coefficients in a certain unitary representation) to the possibility of improving the asymptotic homomorphism. The use of cohomological obstructions goes in essence back to the pioneering work  Kazhdan \cite{art:kazh} on stability of (uniform) approximate representations of amenable groups. The main result of \cite{lubchigletho} is the following theorem.

\begin{theorem}[\cite{lubchigletho}] \label{main2}
Let $\Gamma$ be a finitely presented group such that $$H^2(\Gamma,\mathcal H_\pi)=0$$ for every unitary representation $\pi \colon \Gamma \to {\rm U}(\mathcal H_\pi)$. Then, any asymptotic homomorphism $\varphi_n \colon \Gamma \to {\rm U}(n)$ w.r.t.\ the Frobenius norm is asymptotically close to a sequence of homomorphisms, i.e. $\Gamma$ is Frobenius-stable.
\end{theorem}

It is well-known that a discrete group has Kazhdan's property (T) if and only if $H^1(\Gamma,\mathcal H_{\pi})=0$ for all unitary representations. In \cite{lubchigletho}, the notion of a group to be $n$-Kazhdan was introduced as a vanishing condition of cohomology in dimension $n$ with coefficients in arbitrary unitary representations, see Definition \ref{defvan}. Important work by Garland \cite{MR0320180} and Ballmann-\'Swi\polhk atkowski \cite{MR1465598} provides first examples of $2$-Kazhdan groups. However, those groups all act on Bruhat-Tits buildings of higher rank and thus are residually finite. The remaining delicate work was then to show that nevertheless there \emph{do} exist finitely presented groups which are $2$-Kazhdan and are \emph{not} residually finite. The method in \cite{lubchigletho} is based on Deligne's construction \cite{MR507760} of a non-residually finite central extension of a ${\rm Sp}(2n,\mathbb Z)$. 

Before we outline the definition of the cohomological obstruction to the possibility of improving an asymptotic homomorphism and consider a few examples, let us mention a few open questions.

\begin{question}
Are all amenable (or even all nilpotent or solvable) groups Frobenius-approximated? 
\end{question}

\begin{question}
Is the class of Frobenius-approximated groups closed under central quotients or crossed products by $\mathbb Z$, compare with \cite{MR3320894, MR2566306}?
\end{question}

The analogue of Theorem \ref{main2} also holds for approximation in the sense of (2) above. However, the corresponding cohomology vanishing results in order to apply the theorem in a non-trivial situation are not available. Note that Kirchberg's conjecture discussed above implies that MF-stable groups should be Ramanujan in the sense of \cite{lubsh}.

\subsection{Cohomological obstructions to stability}

In this secion, we want to outline how a cohomological obstruction to stability can be obtained. Consider the family of matrix algebras ${\rm M}_{n}(\mathbb C)$ equipped with some unitarily invariant, submultiplicative norms $\lVert\cdot\rVert_n$, say the Frobenius norms. We consider the ultraproduct Banach space 
\[\Matuf:=\prod_{n\to\uf}( {\rm M}_{n}(\mathbb C),\lVert\cdot\rVert_n),\]
and the metric ultraproduct
\[\Uniuf:=\prod_{n\to\uf}({\rm U}(n),d_{\lVert\cdot\rVert_n}).\]

We can associate an element $[\alpha]\in H^2\left(\Gamma,\prod_{n\to\uf}( {\rm M}_{n}(\mathbb C),\lVert\cdot\rVert) \right)$ to an asymptotic representation $\varphi_n\colon X\to{\rm U}(n)$. This is done so that if $[\alpha]=0$, then the defect can be diminished in the sense that there is an equivalent asymptotic representation $\varphi_n'$ with effectively better defect, more precisely $\defect(\varphi_n')=o_\uf(\defect(\varphi_n))$.

Note that an asymptotic representation as above induces a homomorphism $\varphi_\uf\colon\Gamma\to \Uniuf$ on the level of the group $\Gamma$. Thus $\Gamma$ acts on $\Matuf$ through $\varphi_\uf$. We consider a section $\sigma\colon \Gamma\to \freegrp_k$ of the natural surjection $\freegrp_k \to \Gamma$ and have $\sigma(g)\sigma(h)\sigma(gh)^{-1}\in\langle \! \langle R\rangle \! \rangle$ for all $g,h\in \Gamma$. We set $\tilde\varphi_n = \varphi_n \circ \sigma$. 

Let us now outline how to  define an element in $H^2(\Gamma,\Matuf)$ associated to $\varphi_n$.
To this end we define $c_n:=c_n(\varphi_n):\Gamma\times \Gamma\to {\rm M}_{n}(\mathbb C)$ by
\[c_n(g,h)=
\frac{\tilde\varphi_n(g)\tilde\varphi_n(h)-\tilde\varphi_n(gh)}{\defect(\varphi_n)},\]
for all $n\in\N$ such that $\defect(\varphi_n)>0$ and $c_n(g,h)=0$ otherwise, for all $g,h\in \Gamma$.

Then, it follows that for every $g,h \in \Gamma$, $c_n(g,h)$ is a bounded sequence, so that the sequence defines a map
\[c=(c_n)_{n\in\N}\colon \Gamma\times \Gamma\to\Matuf.\]
The map $c$ is a Hochschild 2-cocycle with values in the $\Gamma$-module $\Matuf$ and $\alpha(g,h):=c(g,h)\varphi_\uf(gh)^*$ is a 2-cocycle in the usual group cohomology.
We call $\alpha$ the cocycle associated to the sequence $(\varphi_n)_{n \in \mathbb N}$.

Assume now that $\alpha$ represents the trivial cohomology class in $H^2(\Gamma,\Matuf)$, i.e.\ there exists a map $\beta\colon\Gamma\to\Matuf$ satisfying
\[\alpha(g,h)=\varphi_\uf(g)\beta(h)\varphi_\uf(g)^*-\beta(gh)+\beta(g),\qquad g,h\in\Gamma.\]
Then, we have $\beta(1_\Gamma)=0$, $\beta(g)=-\varphi_\uf(g)\beta(g^{-1})\varphi_\uf(g)^*$ and $$c(g,h)=\varphi_\uf(g)\beta(h)\varphi_\uf(h)-\beta(gh)\varphi_\uf(gh)+\beta(g)\varphi_\uf(gh).$$
Furthermore, we can choose
$\beta(g)$ to be skew-symmetric for all $g\in\Gamma$.
Now let $\beta$ be as above and let $\beta_n\colon \Gamma\to {\rm M}_{n}(\mathbb C)$ be any bounded and skew-symmetric lift of $\beta$.
Then $\exp(-\defect(\varphi_n)\beta_n(g))$ is a unitary for every $g\in \Gamma$, so we can define a sequence of maps $\psi_n:\Gamma\to{\rm U}(n)$ by 
\[\psi_n(g)=\exp(-\defect(\varphi_n)\beta_n(g))\tilde\varphi_n(g).\]
Note that since $\tilde\varphi_n(1_\Gamma)=1_{n}$ and $\beta_n(1_\Gamma)=0$, we have $\psi_n(1_\Gamma)=1_{n}$.
It follows easily that $\psi_n|_X$ is an asymptotic representation with $$\defect(\psi_n|_X)=O_\uf(\defect(\varphi_n)),$$ but we prove that the defect of $\psi_n|_X$ is actually $o_\uf(\defect(\varphi_n))$. If we define the asymptotic representation $\varphi'_n\colon  X \to{\rm U}(n)$ by $\varphi'_n=\psi_n|_{X}$, the conclusion can be summarized as follows:

\begin{theorem}[\cite{lubchigletho}]\label{speed}
Let $\Gamma=\langle X| R\rangle$ be a finitely presented group and let $\varphi_n\colon X \to{\rm U}(n)$ be an asymptotic representation with respect to a family of submultiplicative, unitarily invariant norms.
Assume that the associated 2-cocycle $\alpha=\alpha(\varphi_n)$ is trivial in $H^2(\Gamma,\Matuf)$. Then there exists an asymptotic representation $\varphi_n'\colon X \to{\rm U}(n)$ such that
$$\dist(\varphi_n,\varphi'_n)=O_\uf(\defect(\varphi_n)) \quad \mbox{and} \quad \defect(\varphi'_n)=o_\uf(\defect(\varphi_n)).$$
\end{theorem}

The converse of Theorem \ref{speed} is also valid in the following sense.
\begin{proposition}[\cite{lubchigletho}]\label{omvendtverden}
Let $\Gamma=\langle X| R\rangle$ be a finitely presented group,
let $\varphi_n,\psi_n: X\to {\mathrm U}(n)$ be asymptotic representations with respect to some family of submultiplicative, unitarily invariant norms and suppose
\begin{itemize}
\item $\dist(\varphi_n,\psi_n)=O_\uf(\defect(\varphi_n))$
and
\item $\defect(\psi_n)=o_\uf(\defect(\varphi_n))$.
\end{itemize}
Then, the 2-cocycle $\alpha$ associated with $(\varphi_n)_{n \in \mathbb N}$ is trivial in $H^2(\Gamma,\Matuf)$.
In particular, if $\varphi_n$ is sufficiently close to a homomorphism for $n$ large enough, then $\alpha$ is trivial.
\end{proposition}

It remains to observe that in case of the Frobenius-norm, the ultraproduct $M_{\mathcal U}$ is a Hilbert space and the action of $\Gamma$ is given by a unitary representation. Together with a somewhat subtle minimality argument this proves Theorem \ref{main2}.

It is now clear that we are in need of large classes of groups for which general vanishing results for the second cohomology with Hilbert space coefficients can be proven. We will discuss some aspects of this problem in the next section.

\subsection{Cohomology vanishing and examples of $n$-Kazhdan groups} \label{hika}
Recall that if $\Gamma$ is a finitely generated group, then $\Gamma$ has Kazhdan's Property (T) if and only if the first cohomology $H^1(\Gamma,\mathcal H_\pi)$ vanishes for every unitary representation $\pi\colon \Gamma\to\Uni(\mathcal H_\pi)$ on a Hilbert space $\mathcal H_\pi$, see \cite{MR2415834} for a proof and more background information.
We will consider groups for which the higher cohomology groups vanish. Higher dimensional vanishing phenomena have been studied in various articles, see for example \cites{MR3284391,  MR1465598
, MR1946553, MR3343347}.

In \cite{lubchigletho}, we proposed the following terminology.
\begin{definition} \label{defvan}
Let $n\in\N$. A group $\Gamma$ is called \emph{$n$-Kazhdan} if $H^n(\Gamma,\mathcal H_\pi)=0$ for all unitary representations $(\pi,\mathcal H_\pi)$ of $\Gamma$. We call $\Gamma$ \emph{strongly $n$-Kazhdan}, if $\Gamma$ is $k$-Kazhdan for $k=1,\ldots, n$.
\end{definition}
So 1-Kazhdan is Kazhdan's classical property (T).
See \cite{MR3284391, MR3343347} for discussions of other related higher dimensional analogues of Property (T). Let's discuss briefly one source of $n$-Kazhdan groups for $n \geq 2$.
Let $K$ be a non-archimedean local field of residue class $q$, i.e.\ if $\mathcal O \subset K$ is the ring of integers and $\mathfrak m \subset \mathcal O$ is its unique maximal ideal, then $q= |\mathcal O/\mathfrak m|$. Let $\mathbf{G}$ be a simple $K$-algebraic group of $K$-rank $r$ and assume that $r \geq 1$. The group $G:= \mathbf G(K)$ acts on the associated Bruhat-Tits building $\mathcal B$. The latter is an infinite, contractible, pure simplicial complex of dimension $r$, on which $G$ acts transitively on the chambers, i.e.\ the top-dimensional simplices. Let $\Gamma$ be a uniform lattice in $G$, i.e.\ a discrete cocompact subgroup of $G$. When $\Gamma$ is also torsion free, then the quotient $X:= \Gamma \backslash \mathcal B$ is a finite $r$-dimensional simplicial complex and $\Gamma= \pi_1(X).$ In particular, the group $\Gamma$ is finitely presented. We will use the following theorem which essentially appears in work of Ballmann-\'Swi\polhk atkowski \cite{MR1465598} building on previous work of Garland \cite{MR0320180}.

\begin{theorem}[Garland, Ballmann-\'Swi\polhk atkowski] \label{higherkazhdan}
For every natural number $r \geq 2$, there exists $q_0(r) \in \mathbb N$ such that the following holds. If $q \geq q_0(r)$ and $G$ and $\Gamma$ are as above, then $\Gamma$ is strongly $(r-1)$-Kazhdan. In particular, if $r\geq 3$, then $\Gamma$ is $2$-Kazhdan.
\end{theorem}

It is very natural to wonder what happens in the analogous real case. It is worth noting that already $H^5({\rm SL}_n(\mathbb Z),\mathbb R)$ is non-trivial for $n$ large enough \cite{borel}; thus ${\rm SL}_n(\mathbb Z)$ fails to be $5$-Kazhdan for $n$ large enough. 
Similarly, note that $H^2({\rm Sp}(2n,\mathbb Z),\mathbb R)=\mathbb R$ for all $n \geq 2$  \cite{borel}, so that the natural generalization to higher rank lattices in real Lie groups has to be formulated carefully; maybe just by excluding an explicit list of finite-dimensional unitary representations.

\begin{question}[\cite{lubchigletho}]
Is it true that ${\rm SL}_n(\mathbb Z)$ is $2$-Kazhdan for $n\geq 4$?
\end{question}

It is worthwhile to return to the remark that Theorem \ref{speed} and the analogue of Theorem \ref{main2} is valid if one replaces $\frob{\cdot}$ with any submultiplicative norm $\lVert\cdot\rVert$ and the assumption that $H^2(\Gamma,V)=\{0\}$ whenever $V=\prod_{n\to\uf}( {\rm M}_{n}(\mathbb C),\lVert\cdot\rVert)$ equipped with some action of $\Gamma$. This, for instance, gives a sufficient condition for stability with respect to the operator norm, but it seems difficult to prove the existence of a group $\Gamma$ with vanishing second cohomology in this case. The following question seems more approachable.

\begin{question} \label{schatt}
Can the above strategy be used to prove stability results w.r.t.\ to the Schatten-$p$-norm?
\end{question}

The techniques rely on submultiplicativity of the norm and thus cannot be directly applied to the normalized Hilbert-Schmidt norm $\hsnorm\cdot$. However, it is worth noting, that since $\frac{1}{\sqrt{k}}\frob{A}=\hsnorm A\leq \opnorm{A}\leq \frob{A}$ for $A\in {\rm M}_k(\mathbb C)$, we get the following immediate corollary to Theorem \ref{main2}.

\begin{corollary}[\cite{lubchigletho}]\label{kazh2hyper}
Let $\Gamma=\langle X | R\rangle$ be a finitely presented 2-Kazhdan group and let $\varphi_n\colon X \to{\rm U}(n)$ be a sequence of maps such that
\[\defect(\varphi_n)=o_\uf(n^{-1/2}),\]
where the defect is measured with respect to either $\hsnorm\cdot$ or $\opnorm\cdot$. Then
$\varphi_n$ is equivalent to a sequence of homomorphisms.
\end{corollary}

The preceding corollary provides some quantitative information on Connes' Embedding Problem. Indeed, if a finitely presented, non-residually finite, 2-Kazhdan group is Connes-embeddable, then there is some upper bound on the quality of the approximation in terms of the dimension of the unitary group. Needless to say it would be very interesting to decide if groups as above are Connes-embeddable. A positive answer to Question \ref{schatt} for $p> 2$ should lead to improvements in Corollary \ref{kazh2hyper}.

\section{Applications to group theory}

\subsection{The basic setup}

For any group $\Gamma$, an element $w$ in the free product $\Gamma \ast {\mathbf F}_k$  determines a word map $w \colon  \Gamma^{\times n} \to \Gamma$ given by evaluation. Let us denote by $\varepsilon \colon  \Gamma \ast {\mathbf F}_k \to {\mathbf F}_k$ the natural augmentation which sends $\Gamma$ to the neutral element and call $\varepsilon(w)$ the {\it content} of $w$. We call $w \in \Gamma \ast {\mathbf F}_k$ a group word in $k$ variables with coefficients in $\Gamma$. Every group word $w \in \Gamma \ast {\mathbf F}_k$ determines an equation $w(X)=1$ in $k$ variables with coefficients in $\Gamma$ in an obvious way. We say that $w \in \Gamma \ast {\mathbf F}_k$ can be solved {\it over} $\Gamma$ if there exists an overgroup $\Lambda \supseteq \Gamma $ and $g_1,\dots,g_k \in \Lambda$ such that $w(g_1,\dots,g_k)=1$, where $1$ denotes the neutral element in $\Lambda$. Similarly, we say that it can be solved {\it in} $\Gamma$ if we can take $\Lambda=\Gamma$. It is clear that an equation $w \in \Gamma \ast {\mathbf F}_k$ can be solved over $\Gamma$ if and only if the natural homomorphism $\Gamma \to \Gamma \ast {\mathbf F}_k/\langle\!\langle w \rangle\! \rangle$ is injective. Similarly, an equation can be solved in $\Gamma$ if and only if the natural homomorphism $\Gamma \to \Gamma \ast {\mathbf F}_k/\langle\!\langle w \rangle\! \rangle$ is split-injective, i.e., it has a left inverse. 

The study of equations over groups dates back to the work of Bernhard Neumann \cite{MR0008808}. There is an extensive literature about equations over groups, including \cites{ MR919828, MR0166296, MR614523, MR1218513, MR0142643, MR0008808, MR3604379, MR3043434}. In this section, we plan to survey some observations and results that were obtained in joint work with Klyachko \cite{MR3604379}.

Let us start with an observation. It is well-known that not all equations with coefficients in $\Gamma$ are solvable over $\Gamma$. For example if $\Gamma=\langle a,b | a^2,b^3 \rangle$, then the equation $w(x) = xax^{-1}b$ with variable $x$ is not solvable over $\Gamma$. Indeed, $a$ and $b$ cannot become conjugate in any overgroup of $\Gamma$. Another example involving only one kind of torsion is $\Gamma = \Z/p\Z = \langle a \rangle$ with the equation $w(x) = xax^{-1} a xa^{-1}x^{-1}a^{-2}$. 
However, in both cases we have $\varepsilon(w)=1 \in {\mathbf F}_k$. Indeed, the only known examples of equations which are not solvable over some $\Gamma$ are equations whose content is trivial. 
We call an equation $w \in \Gamma \ast {\mathbf F}_k$ singular if its content is trivial, and non-singular otherwise. This leads to the following question:

\begin{question}[\cite{MR3604379}] \label{conj}
Let $\Gamma$ be a group and $w \in \Gamma \ast {\mathbf F}_k$ be an equation in $n$ variables with coefficients in $\Gamma$. If $w$ is non-singular, is it true, that it is solvable over $\Gamma$? In addition, if $\Gamma$ is finite, can we find a solution in a finite extension?
\end{question}

The case $k=1$ is the famous Kervaire-Laudenbach Conjecture.
The one-variable case was studied in work by Gerstenhaber-Rothaus, see \cite{MR0166296}. They showed that if $\Gamma$ is finite, then every non-singular equation in one variable can be solved in a finite extension of $\Gamma$. Their proof used computations in cohomology of the compact Lie groups ${\rm U}(n)$. Their strategy was to use homotopy theory to say that the associated word map $w \colon {\rm U}(n) \to {\rm U}(n)$ has a non-vanishing degree (as a map of oriented manifolds) and thus must be surjective. Any preimage of the neutral element provides a solution to the equation $w$. The key to the computation of the degree is to observe that the degree depends only on the homotopy class of $w$ and thus -- since ${\rm U}(n)$ is connected -- does not change if $w$ is replaced by $\varepsilon(w)$. The computation of the degree is now an easy consequence of classical computations of Hopf \cite{hopf}. We conclude that any non-singular equation  in one variable with coefficients in ${\rm U}(n)$ can be solved ${\rm U}(n)$ -- thus ${\rm U}(n)$ deserves to be called \emph{algebraically closed}.

The property of being algebraically closed is easily seen to pass to arbitrary Cartesian products of groups and arbitrary quotients of groups. As a consequence, non-singular equations in one variable with coefficients in $\Gamma$ as above can be solved over $\Gamma$ if $\Gamma$ is isomorphic to a subgroup of a quotient of the infinite product $\prod_n \! {\rm U}(n)$ -- an observation that is due to Pestov \cite{MR2460675}:

\begin{theorem}[\cite{MR2460675}]
The Kervaire-Laudenbach Conjecture holds for Connes-embeddable groups.
\end{theorem}

Note that this covers all amenable groups, or more generally, all sofic groups \cite{MR2460675}. As we have discussed, the Connes' Embedding Conjecture predicts (among other things) that every countable group is Connes-embeddable and thus implies the Kervaire-Laudenbach Conjecture -- this was also observed by Pestov in \cite{MR2460675}.

\vspace{0.2cm} Actually, Gerstenhaber-Rothaus \cite{MR0166296} studied the more involved question whether $m$ equations of the form $w_1,\dots,w_m \in \Gamma \ast {\mathbf F}_k $ in $k$ variables can be solved simultaneously over $\Gamma$. Their main result is that this is the case if $\Gamma$ is finite (or more generally, locally residually finite) and the presentation two-complex $X := K\langle X | \varepsilon(w_1),\dots,\varepsilon(w_m) \rangle$ satisfies $H_2(X,\mathbb Z)=0$, i.e., the second homology of $X$ with integral coefficients vanishes. 

Later, Howie \cite{MR614523} proved the same result for locally indicable groups and conjectured it to hold for all groups -- we call that Howie's Conjecture. Again, Connes' Embedding Conjecture implies Howie's Conjecture -- and more specifically, every Connes-embeddable group satisfies Howie's Conjecture.

\subsection{Topological methods to prove existence of solutions}

The main goal of \cite{MR3604379} was to provide examples of singular equations in many variables which are solvable over every Connes-embeddable group, where the condition on the equation {\it only} depends on its content. Indeed, we gave a positive answer to Question \ref{conj} when $k =2$ in particular cases. This should be compared for example with results of Gersten \cite{MR919828}, where the conditions on $w$ depended on the unreduced word obtained by deleting the coefficients from $w$. 

\begin{theorem}[\cite{MR3604379}] \label{mainthm}
Let $\Gamma$ be a Connes-embeddable group. An equation in two variables with coefficients in $\Gamma$ can be solved over $\Gamma$ if its content does not lie in $[\F_2,[\F_2,\F_2]]$. Moreover, if $\Gamma$ is finite, then a solution can be found in a finite extension of $\Gamma$.
\end{theorem}

In order to prove our main result we had to refine the study of Gerstenhaber-Rothaus on the effect of word maps on cohomology of compact Lie groups. 
Again, the strategy is to show that such equations can be solved in $\SU(n)$ for sufficiently many $n \in \N$. More specifically, we proved:

\begin{theorem}[\cite{MR3604379}] \label{sup}
Let $p$ be a prime number and let $w \in \SU(p) \ast {\mathbf F}_2$ be a group word. If $$\varepsilon(w) \not \in [{\mathbf F}_2,{\mathbf F}_2]^p [{\mathbf F}_2,[{\mathbf F}_2,{\mathbf F}_2]],$$ then the equation $w(a,b)=1$ can be solved in $\SU(p)$.
\end{theorem}

If $\varepsilon(w) \not \in [\F_2,\F_2]$, then this theorem is a direct consequence of the work of Gerstenhaber-Rothaus. However, if $\varepsilon(w) \in [\F_2,\F_2]$, then a new idea is needed. We showed -- under the conditions on $p$ which are mentioned above -- that the induced word map $w \colon  \PU(p) \times \PU(p) \to \SU(p)$ is surjective, where $\SU(p)$ denotes the special unitary group and $\PU(p)$ its quotient by the center. Again, the strategy was to replace $w$ by the much simpler and homotopic map induced by $\varepsilon(w)$ and study its effect on cohomology directly. 

The proof is a \emph{tour de force} in computing the effect of $\varepsilon(w) \colon \PU(p) \times \PU(p) \to \SU(p)$ in cohomology with coefficients in $\mathbb Z/p\mathbb Z$. Using fundamental results of Serre, Bott, and Baum-Browder on the $p$-local homotopy type of spheres, lens spaces and projective unitary groups and  finally computations of Kishimoto and Kono \cite{MR2544124}, we managed to show that the image of the top-dimensional cohomology class of $\SU(p)$ is non-trivial. This implies that no map that is homotopic to $\varepsilon(w)$ can be non-surjective. In particular, we can conclude as in the arguments of Gerstenhaber-Rothaus that $w\colon \PU(p) \times \PU(p) \to \SU(p)$ must be surjective.

In general, the assumption on $\varepsilon(w)$ cannot be omitted in the previous theorem. Indeed, in previous work the following result (independently obtained by Elon Lindenstrauss) was shown.

\begin{theorem}[\cite{MR3283561}]
For every $k \in \N$ and every $\varepsilon>0$, there exists $w \in {\mathbf F}_2 \setminus \{e\}$, such that 
$$\|w(a,b)-1_n\| < \varepsilon,\quad \forall a,b \in {\rm SU}(k).$$
In particular, the equation $w(a,b)=g$ is not solvable when $\|g-1_n\|\geq \varepsilon$.
\end{theorem}

The construction that proves the preceding theorem yields words in $\F_2$ that lie deep in the derived series, so that there is no contradiction with Theorem \ref{sup}. 

The surjectivity of word maps without coefficients is an interesting subject in itself. Larsen conjectured that for each non-trivial $w \in \F_2$ and $n$ high enough, the associated word map $w \colon  \PU(n) \times \PU(n) \to \PU(n)$ is surjective. This was shown (with some divisibility restrictions on $n$) for words not in the second derived subgroup of $\F_2$ by Elkasapy and the author in \cite{thomelk}.

In a similar direction, we believe that for $n$ high enough -- or again, with some divisibility restrictions -- the word map $w$ should define a non-trivial homotopy class and be not even homotopic to a non-surjective map.
In this context let us mention some questions that appear naturally at the interface between homotopy theory and the study of word maps.
Given a topological group $G$, it is natural to study the group of words modulo those which are null-homotopic. Indeed, we set $$N_{n,G} := \left\{w \in {\mathbf F}_k | w \colon G^n \to G \mbox{ is
homotopically trivial} \right\}$$ and define ${\mathcal H}_{n,G} := {\mathbf F}_k/N_{n,G}$.
\begin{question}
Can we compute ${\mathcal H}_{2,{\rm SU}(n)}$?
\end{question}
See \cite{MR0106465, MR1233412} for partial information about  $\mathcal H_{n,G}$ in special cases. For example, it follows from classical results of Whitehead that ${\mathcal H}_G$ is $k$-step nilpotent for some $k \leq 2 \dim(G)$.

Similarly, we call $w \in {\mathbf F}_k$ {\it homotopically surjective} with respect to $G$ if every map in the homotopy class of $w \colon  G^{\times n} \to G$ is surjective. 

\begin{question} Let $w \in {\mathbf F}_2$ be non-trivial.
Is $w \colon \PU(n) \times \PU(n) \to \PU(n)$ homotopically surjective for large $n$?
\end{question}

In order to study words which lie deeper in the lower central series, we suspected in \cite{MR3604379} that it might be helpful to oberserve that the induced word map $w \colon  \PU(p) \times \PU(p) \to \PU(p)$ does not only lift to $\SU(p)$  -- which is the simply connected cover of $\PU(p)$ -- but lifts even to higher connected covers of $\PU(p)$. Indeed, for example one can show that if $w \in [\F_2,[\F_2,\F_2]]$, then the associated word map lifts to the complex analogue of the string group.

\section*{Acknowledgments}
This research was supported by ERC Consolidator Grant No.\ 681207. I thank Jakob Schneider for careful reading.

\end{document}